\documentclass[11pt]{article}
\pdfoutput=1
\usepackage[ascii]{inputenc}
\usepackage[T1]{fontenc}
\usepackage{amsmath,amssymb,amsfonts,textcomp,epsfig,amsthm}
\usepackage{color}
\usepackage{calc}
\usepackage[all]{xy}
\usepackage{graphicx}
\usepackage{tikz}

\newcommand{\R}{\mathbb{R}}

\newtheorem{thm}{Theorem}

\newtheorem{defi}{Definition}

\title{{\bf Centers of reversible cubic  perturbations of the  symmetric 8-loop Hamiltonian  }}
\author{
Fatma Sassi\\
Facult\'e des Sciences de Sfax, D\'epartement de Math\'ematiques,\\
BP 1171, 3000 Sfax, Tunisie.\\
E-mail: \texttt{fatma.lazher@gmail.com}\\
\\
Ameni Gargouri\\
Facult\'e des Sciences de Sfax, D\'epartement de Math\'ematiques,\\
BP 1171, 3000 Sfax, Tunisie.\\
E-mail: \texttt{ameni.gargouri@gmail.com}\\
\\
Lubomir Gavrilov\\
Institut de Math\'{e}matiques de Toulouse, UMR 5219\\
Universit\'{e}  de Toulouse,  31062 Toulouse,  France.\\
E-mail: \texttt{lubomir.gavrilov@math.univ-toulouse.fr}\\
\\
Bassem Ben Hamed\\
Ecole Nationale d'Electronique et des T\'el\'ecommunication de Sfax, \\
Route de Tunis, Cit\'e ElOns, Technop\^ole de Sfax - 3018 Sfax, Tunisie.\\
E-mail: \texttt{bassem.benhamed@enetcom.usf.tn}}

\date{}
\begin{document}

\maketitle

\begin{abstract}
We show that the center set  of reversible cubic systems, close to the symmetric Hamiltonian system $x'=y, y'= x-x^3$ has two irreducible components of co-dimension two in the parameter space. One of them corresponds to the Hamiltonian stratum, the other to systems which are polynomial  pull back of an appropriate linear system
 \\  
{\it Keywords:} center-focus problem, cubic vector field.  
\end{abstract}

\section{Introduction}
\label{section1}
\begin{figure}\label{Fig1}
\begin{center}
\includegraphics[width=8cm]{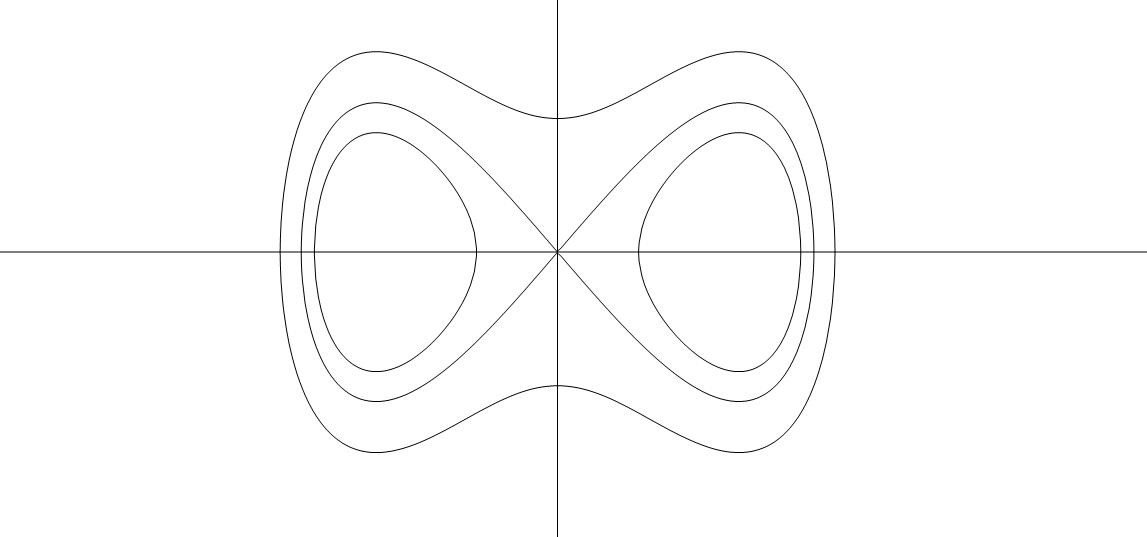}
\end{center}
\caption{The phase portrait of  (\ref{nonperturbed}).}
\end{figure}

The paper is a contribution to the study of the center set of plane polynomial cubic vector fields, in a neighbourhood of the symmetric cubic vector field, see Fig.\ref{Fig1},
\begin{eqnarray}\label{nonperturbed}
X_{0}:
\left\{\begin{array}{ccl} \dot{x}&=&y \\
\dot{y}&=& x-x^{3}
\end{array}\right.
\end{eqnarray}
$X_0$ is Hamiltonian, and has a first integral 
\begin{equation*}
H(x,y)= \frac{1}{2}y^{2}- \frac{1}{2}x^{2} +\frac{1}{4}x^{4}.
\end{equation*}
Under a small analytic perturbation the centres of $X_0$ near $(1,0)$ and $(-1,0)$ are either simultaneously destroyed, or simultaneously persistent. The set of cubic  vector fields close to $X_0$ and having 
 a center near $(\pm1,0)$  is 
 the \emph{center set} $\mathcal C_0$. It is known that the center set is an algebraic set in the space of parameters, but even the number of its irreducible components is unknown.

Recently Iliev, Li and Yu  \cite{Iliev}  studied special one-parameter families of perturbations of the form
\begin{eqnarray*}
X_{\varepsilon}:
\left\{\begin{array}{ccl} \dot{x}&=&y + \varepsilon P(x,y) \\
\dot{y}&=& x-x^{3} + \varepsilon Q(x,y)
\end{array}\right.
\end{eqnarray*}
where $P,Q$ are arbitrary \emph{fixed} real cubic polynomials. The displacement function near the singular points $(\pm1,0)$ has an analytic expansion
\begin{align}
d(h,\varepsilon)= \varepsilon M_1(h)+ \varepsilon^2 M_2(h) + \varepsilon^3 M_3(h) + \dots
\end{align}
where as usual $h$ is the restriction of the Hamiltonian $H$ on a cross-section to the vector field. The so called Melnikov functions $M_k$ vanish if and only if the displacement map is the zero map, that is to say the the centers  $(1,0)$ and $(-1,0)$ are simultaneously persistent.
 It was shown then in \cite[Theorem 1]{Iliev} that if $M_1=M_2=M_3=M_4=0$, then the displacement map is identically zero and therefore $X_\varepsilon$ has a center near $(1,0)$ and $(-1,0)$. The set of such system is an algebraic set $\mathcal C^l$ contained in $\mathcal C$. It turns out that $\mathcal C^l$ is a union of vector spaces, and a vector field $X_\varepsilon$ which belongs to an irreducible component of  $\mathcal C^l$ is either Hamiltonian, or $y$-reversible, or $x$-reversible.
It is clear that when a vector field $X_\varepsilon$  is Hamiltonian, or $y$-reversible, then it has a center near $(\pm1,0)$. If, however, the vector field is $x$-reversible, it does not follow that it has a center near $(\pm1,0)$. Therefore, it makes a sense to consider the case of $x$-reversible systems separately. 

The purpose of this paper is to give complete description of the center set $\mathcal C$ under the restriction that the vector field is $x$-reversible, that is to say, the associated foliation by orbits is invariant under the involution $(x,y)\to (-x,y)$. Our approach is the following. The invariance under $x \to -x$ suggests to introduce the quotient vector field which is stil polynomial and cubic. It turns out, that this new vector field is of Li\'enard type, whose centers were extensively studied. We apply a classical result of \cite[Cherkas] {cherkas} revisited recently by \cite[Chrystopher]{chri99}. As a result we obtain that in this case the center set has two co-dimension two smooth irreducible components which correspond either to Hamiltonian systems, or to systems obtained as polynomial pull back from linear systems. From this, the result of the paper follows.


\begin{figure}\label{Fig2}
\begin{center}
\includegraphics[width=10cm]{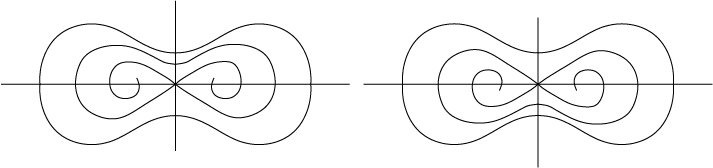}
\end{center}
\caption{The phase portrait of the perturbed reversible system \eqref{reversible} with two foci. The exterior period annulus persists.}
\end{figure}

\section{Statement of the  results }
\label{section2}

Consider the following perturbed cubic system  
\begin{eqnarray}\label{perturbed}
X_{\lambda}:
\left\{\begin{array}{ccl} \dot{x}&=&y + P(x,y) \\
\dot{y}&=& x-x^{3} + Q(x,y)
\end{array}\right.
\end{eqnarray}
where  
$$
P(x,y)= \sum_{i+j\leq 3} a_{ij} x^iy^j, \; Q(x,y)= \sum_{i+j\leq 3} b_{ij} x^iy^j
$$
$\lambda = \{ a_{ij},b_{ks} \}$ are small parameters. For $\lambda=0$ the system $X_0$ has a first integral 
\begin{equation*}
H(x,y)= \frac{1}{2}y^{2}- \frac{1}{2}x^{2} +\frac{1}{4}x^{4} 
\end{equation*}
and two centers at $(x,y)=(\pm1,0)$ shown on Fig 1.
The perturbed vector field $X_\lambda$ has therefore a saddle, close to the origine $(0,0)$, as well two anti-saddles (centers or foci) close to $(\pm1,0)$. 
The anti-saddles near $(\pm1,0)$ are either simultaneously centers, or saddles\footnote{this fact is obvious if the vector field $X_\lambda$ is $x$-reversible, but holds true for arbitrary analytic perturbations too}.

The center set  $\mathcal C_0$ of \emph{small} parameters $\lambda$ for which the vector filed has a center near $(\pm1,0)$ is a germ of analytic set in the space of parameters $\lambda$.
The set $\mathcal C_0$ is in fact algebraic, it is globally defined as the zero set of a finite family of polynomials in $\lambda$, but its number of irreducible components is not known in general. 

The purpose of the present paper is to describe $\mathcal C_0$ in the particular case, when $X_\lambda$ is reversible in $x$. More precisely
\begin{defi}
The vector field (\ref{perturbed}) is said to be reversible with respect to $x$ provided that the involution $x\mapsto -x$ sends $X_\lambda$ to $-X_\lambda$, or equivalently
 $P(-x,y)=P(x,y) $, $Q(-x,y)=-Q(x,y)$.  The set of reversible cubic systems (\ref{reversible})  having a center near $(\pm 1, 0)$ is denoted  $\mathcal C^R_0$.
\end{defi}
The set of reversible in $x$ vector fields $X_\lambda^R$ 
\begin{eqnarray}\label{reversible}
X_{\lambda}^R:
\left\{\begin{array}{ccl} \dot{x}&=&y + P^R(x,y) \\
\dot{y}&=& x-x^{3} + Q^R(x,y)
\end{array}\right.
\end{eqnarray}
is therefore parameterised by the space of special cubic polynomials $P^R,Q^R$ of the form
\begin{align}
\label{p}
P^R(x,y)  &= a_{00}+a_{20}x^2+a_{21}x^{2}y +a_{01}y+a_{02}y^{2}+ a_{03}y^{3}, \\ 
\label{q}
Q^R(x,y)&= b_{10}x+b_{11}xy+b_{12}xy^{2}+b_{30}x^{3}
\end{align}
where $\lambda = \{ a_{ij},b_{ks} \}$ are small real parameters.
 The possible phase portraits of  
the perturbed vector field $X_\lambda^R$ in the finite plane (that is to say, in a disc of finite radius) is shown on Fig.1, where there is some unknown number of limit cycles. We note the orbits which belong to  the exterior period annulus are always closed, because the vector field $X_\lambda^R$ in consideration is reversible.  We have a strict inclusion $\mathcal C^R_0 \subset \mathcal C_0$ and $\mathcal C^R_0$ is  an algebraic set in the parameter space.
We shall prove

\begin{thm}
\label{th0}
The center set $\mathcal C^R_0$  of reversible cubic systems $X_\lambda^R$ with a center near $(\pm1,0)$ has two  irreducible components of co-dimension two in the set of polynomials  (\ref{p}), (\ref{q}). The components correspond either to Hamiltonian systems or  to systems which are obtained as polynomial  pull back from an appropriate linear system.
\end{thm}
To describe \emph{explicitly }  $\mathcal C^R_0$  we shall normalise first $X_\lambda$ as follows. Note that the affine transformations 
$$
(x,y) \mapsto (\alpha x, y + \beta)
$$
transform a reversible cubic system to a reversible cubic system of the same form, and therefore act on the parameter space and the center set  $\mathcal C^R_0$. 

Therefore, performing an appropriate affine change of $x,y$, we may assume that $X_\lambda^R$ has a singular point at $(1,0)$ for all sufficiently small $\lambda$ and by the $x$-reversibility, it will have another singular point at $(-1,0)$.  The normalised vector field (\ref{reversible}) takes the form
\begin{eqnarray}\label{normal}
X_{\lambda}^R:
\left\{\begin{array}{ccl} \dot{x}&=&y +a_{20}(x^{2}-1)+ a_{21}x^{2}y  + a_{01}y + a_{02}y^{2}+ a_{03}y^{3}  \\
\dot{y}&=& x-x^{3} +b_{30}(x^{3} -x) + b_{11}xy +b_{12}xy^{2}
\end{array}\right.
\end{eqnarray}
(we denote the coefficients of this normalised reversible vector field by the same letters $a_{ij}$).
Theorem \ref{th0} is an obvious consequence of the following
\begin{thm} \label{thm1}
The system (\ref{normal}) 
 has a non-degenerate real center at $(\pm1,0) $ if and only if  
\begin{align*} 2a_{20}+b_{11} = 0
\end{align*}
and either
\begin{itemize}
\item $a_{21} + b_{12}=0$ (Hamiltonian case), or
\item  $(1-b_{30})a_{02}=a_{20}(2a_{21}- b_{12})$  (pull back case).
\end{itemize}
In the second case the system  (\ref{normal})  is a polynomial pull back of a linear system under the map
\begin{align*}
(x,y) \to ( x^2-1+P_2(y), y^2), P_{2}(y)=  \frac{1}{b_{30}-1}( b_{11}y +b_{12}y^{2}) .
\end{align*}
\end{thm}
Note that the trace of $X_\lambda^R$ at $(\pm1,0) $ equals $2a_{20}+b_{11}$. 
To determine the center conditions of (\ref{normal})
we shall use the Cherkas-Christopher theorem which we explain in the next section.

\section{The Cherkas-Christopher Theorem}
 \label{section3}
Consider the plane cubic differential system
\begin{eqnarray}\label{cubic1}
\left\{\begin{array}{ccl} \dot{\xi}&=&P_{3}(y)+ P_{1}(y) \xi  \\
\dot{y}&=&  -\xi -P_{2}(y)
\end{array}\right.
\end{eqnarray}
where
\begin{align}
P_{1}(y)= a_{0}+ a_{1}y,\; P_{2}(y)=b_{1}y+ b_{2}y^{2}, \; P_{3}(y)=c_{1}y+ c_{2}y^{2}+ c_{3}y^{3} \in \R[y] .
\end{align}
Upon substituting $\xi \to R(x,y)$ in  (\ref{cubic1}), where $R$ is a quadratic polynomial, we get a new cubic system  ($\ref{cubic1}^*$), which is pull back of (\ref{cubic1}) under the polynomial map $(\xi,y)= (R(x,y),y)$. In such a way centers of (\ref{cubic1}) produce new centers of more general cubic  systems.

In this section we determine the necessary and sufficient conditions so that the cubic system (\ref{cubic1})
 has a non-degenerate singular point at the origin of center type. We assume therefore through this section, that (\ref{cubic1}) has already a linear center at the origin, that is to say
 \begin{align}
 a_0-b_1 = 0, c_1 - a_0b_1 >0 .
 \end{align}
Substituting $\xi  = u-P_{2}(y) $ in (\ref{cubic1}) we  get:
\begin{equation*}
udu+\left[ u \left( (a_{0}-b_{1}) + (a_{1}-2b_{2})y \right) +(c_{1}-a_{0}b_{1})y +(c_{2}-a_{0}b_{2}-a_{1}b_{1})y^{2}+(c_{3}-a_{1}b_{2})y^{3}   \right] dy=0,
\end{equation*}
which is the polynomial foliation of the Li\'enard equation
\begin{eqnarray}\label{cubic2}
\left\{\begin{array}{ccl} \dot{y}&=&u  \\
\dot{u}&=& -q(y)-up(y)
\end{array}\right.
\end{eqnarray}
where
\begin{align}
p(y)&=a_{0}-b_{1} + (a_{1}-2b_{2})y , \;\\ q(y)&= 
(c_{1}-a_{0}b_{1})y +(c_{2}-a_{0}b_{2}-a_{1}b_{1})y^{2}+(c_{3}-a_{1}b_{2})y^{3} .
\end{align}
The center set of the above Li\'enard system is well known, see \cite{lubo} and \cite{cherkas}.
It has a center at the origin if and  only if the primitives
$P(y)= \int p(y) dy $ and $ Q(y)=\int q(y)dy$ have a common composition factor $W(y)$ with a Morse critical point at the origin. Therefore 
$a_{0}-b_{1}=0$, $W=P$ and there exists a degree two polynomial $\tilde Q$, such that $Q(y)= \tilde Q (W(y))$.The latter is equivalent to the condition that 
$Q$ is even, when $p\neq 0$. In the case $p=0$ the system is obviously Hamiltonian. We resume the result in the following
\begin{thm}
\label{th3}
The system of Li\'enard type (\ref{cubic1}) has a non-degenerate center at the origin if and only if
\begin{align*}
a_{0}-b_{1}=0, \; c_1 - a_0 b_1>0 
\end{align*}
and either
\begin{itemize}
\item $a_1-2b_2=0$ (Hamiltonian case), or
\item $c_{2} -a_{0}(a_{1}+b_{2}) =0 $ (pull back case) .
\end{itemize}
\end{thm}

\section{Proof of Theorem \ref{thm1}} 
\label{section4}
The substitution $\xi= x^{2}-1$ takes the system (\ref{normal}) 
\begin{align*}
(y +a_{20}(x^{2}-1)+ a_{21}x^{2}y  + a_{01}y + a_{02}y^{2}+ a_{03}y^{3}) dy  
-(x-x^{3} +b_{30}(x^{3} -x) + b_{11}xy +b_{12}xy^{2})dx = 0
\end{align*}
to the Li\'enard form (\ref{cubic1}) 
\begin{align*}
(\xi +P_{2}(y)) d\xi+ (P_{3}(y)+ P_{1}(y) \xi )dy =0
\end{align*}
where
\begin{align*}
P_{1}(y)&= \frac{2}{1-b_{30}}(a_{20}+ a_{21}y) \\
P_{2}(y)&= \frac{1}{b_{30}-1}( b_{11}y +b_{12}y^{2}) \\
 P_{3}(y) &= \frac{2}{1-b_{30}}((1+a_{01}+a_{21})y+a_{02}y^{2}+a_{03}y^{3})
 \end{align*}
 
Applying  Theorem \ref{th3} with 
\begin{align*}
&a_0= \frac{2a_{20}}{1-b_{30}}, &a_1&=\frac{2a_{21}}{1-b_{30}} &  \\
&b_1= -\frac{b_{11} }{1-b_{30}}, &b_2 &=  -\frac{b_{12} }{1-b_{30}} & \\
&c_1= \frac{2(1+a_{01}+a_{21})}{1-b_{30}} , &c_2&= \frac{2 a_{02} }{1-b_{30}}, &c_3 = \frac{2 a_{03} }{1-b_{30}}
\end{align*}
we obtain the equation for the center set $\mathcal C_0^R$ for the normalised  reversible vector field (\ref{normal}).

Finally, to find the first integral in the logarithmic case, we recall that (\ref{normal}) is a polynomial pull back 
under $ u=x^2-1 +P_2(y)$
of the Li\'enard system
\begin{eqnarray}\label{cubic3}
\left\{\begin{array}{ccl} \dot{y}&=&u  \\
\dot{u}&=& -
q(y)-up(y)
\end{array}\right.
\end{eqnarray}
where (in the pull back case)
\begin{align*}
p(y)&=  (a_{1}-2b_{2})y , \;\\ q(y)&= 
(c_{1}-a_{0}b_{1})y +(c_{3}-a_{1}b_{2})y^{3} .
\end{align*}
As (\ref{cubic3}) is also  $y$-reversible then it is a pull back of the following linear system 
\begin{eqnarray}\label{linear}
\left\{\begin{array}{ccl} \dot{v}&=&2u  \\
\dot{u}&=& -(c_{1}-a_{0}b_{1}) - (c_{3}-a_{1}b_{2})v
-u ( a_{1} -2b_{2})
\end{array}\right.
\end{eqnarray}
under the map $y\to v =y^2$. This completes the proof of Theorem \ref{thm1}.

We note finally that he system (\ref{linear}) (as any non-degenerate linear system) has a logarithmic first integral of the form $l_1^{\alpha}l_2^{\beta}$ where $l_1,l_2$ are linear functions in $\lambda,u$ and $\alpha, \beta$ are suitable complex   numbers. Although (\ref{linear}) has a Darboux type first integral, it has no center, except in the Hamiltonian case $a_{1} -2b_{2}=0$. 
The reversible vector field $X_\lambda^R$ has also a first integral of Darboux type (pull back of the Darboux first integral of the linear system), but its centers near $(\pm1,0)$ are of pull back type. An explicit computation of this integral in some cases can be found in \cite[section 5]{Iliev}.
$\Box$

\end{document}